\documentclass[reqno,11pt]{amsart}

\usepackage{enumerate}
\usepackage{mdwlist}
\newtheorem{thm}{Theorem}[section]
\newtheorem{cor}[thm]{Corollary}
\newtheorem{lem}[thm]{Lemma}
\newtheorem{prop}[thm]{Proposition}

\theoremstyle{definition}

\newtheorem{rem}[thm]{Remark}

\numberwithin{equation}{section}

\def\iso{\textrm{iso}}

\begin{document}

\title{On the spectral $\nu$-continuity}

\author{S. S\'anchez-Perales}
\address{Universidad Tecnol\'ogica de la Mixteca, Instituto de F\'isica y Matem\'aticas, Km. 2.5 Carretera a Acatlima, 69000 Oaxaca, Mexico.}
\email{es21254@yahoo.com.mx}

\author{S. Palafox}
\address{Universidad Tecnol\'ogica de la Mixteca, Instituto de F\'isica y Matem\'aticas, Km. 2.5 Carretera a Acatlima, 69000 Oaxaca, Mexico.}
\email{sergiopalafoxd@gmail.com}

\author{T. P\'erez-Becerra}
\address{Universidad Tecnol\'ogica de la Mixteca, Divisi\'on de Estudios de Postgrado, Km. 2.5 Carretera a Acatlima, 69000 Oaxaca, Mexico}
\email{tompb55@hotmail.com}


\begin{abstract}
In this paper we study the $\nu$-continuity of the spectrum and some of its parts. We show that the approximate point spectrum $\sigma_{ap}$ is upper semi-$\nu$-continuous at every Fredholm operator, then we give sufficient conditions to guarantee the $\nu$-continuity of $\sigma_{ap}$. Also we show that the restriction of the Weyl spectrum on the class of essentially $G_1$ operators is $\nu$-continuous. Finally, we investigate the $\nu$-continuity of the spectrum on the class of $p$-hyponormal operators.
\end{abstract}


\maketitle

\section{Introduction}
The spectral  continuity is a relevant subject in Banach theory and operator theory. Several authors have studied this topic using different types of convergence. Regarding norm convergence, Newburgh in \cite{MR0051441} was one of the first to do  systematically investigated in general Banach algebras. More completely results in the case of the algebra $B(H)$ for $H$ a Hilbert space are given in the series of papers of Conway and Morrel, see \cite{MR543882,MR634170,MR701023}. Other articles on this topic are \cite{MR1212728,MR1457127,MR1814478,MR1793814,MR2894894,sanchez4}. Another type of convergence for the study of spectral continuity is the $\nu$-convergence given by Ahues in \cite{MR1886113}. Also developed by  S\'anchez et al. \cite{MR3388799, MR3626681,MR4043866} and Ammar et al. \cite{MR3624565, MR4075017}.

In this paper we study the $\nu$-continuity of the spectrum function. In section 2, we  review this concept on arbitrary unital Banach algebras. We show that in an unital abelian Banach algebra the spectrum function is $\nu$-continuous. In section 3, we study the $\nu$-continuity of the spectrum and the approximate point spectrum  on the algebra $B(X)$, for a Banach space $X$. We show that for the particular case of a Hilbert space $H$, the $\nu$-continuity of spectrum is equivalent under certain conditions to the usual continuity of the spectrum.   Also, we show that the approximate point spectrum is upper semi-$\nu$-continuous at every Fredholm operator, then we give sufficient conditions to guarantee the $\nu$-continuity of this function. In section 4, we show that the restriction of the Weyl spectrum on the class of essentially $G_1$ operators is $\nu$-continuous. Finally, we investigate the $\nu$-continuity of the spectrum on the class of $p$-hyponormal operators.

The next concepts are part of classical point set topology.  Let  $(E_n)$ be a sequence of arbitrary subsets of $\mathbb{C}$ and define the limits inferior and superior of
$(E_n)$ as follows:
\begin{itemize}
\item[]  $\liminf E_n=\{\lambda\in \mathbb{C}\,|\,$ for every $\epsilon>0$, there exists $N\in\mathbb{N}$ such that $B(\lambda, \epsilon)\cap E_n\neq \emptyset$ for all   $n\geq N \}$.
\item[] $\limsup E_n=\{\lambda\in \mathbb{C}\,|\,$ for every $\epsilon>0$, there exists $J\subseteq \mathbb{N}$  infinite such that $B(\lambda, \epsilon)\cap E_n\neq \emptyset$ for all  $n\in J \}$.
\end{itemize}

\begin{rem}\label{rem2}
Let $(E_n)$ be a sequence of non-empty subsets of $\mathbb{C}$. The following properties are hold:
\begin{enumerate}[(a)]
\item $\liminf E_n$ and $\limsup E_n$ are closed subsets of $\mathbb{C}$.
\item $\lambda\in\limsup E_n$ if and only if there exists an increasing sequence of natural numbers $n_1<n_2<n_3<\cdots$\, and points $\lambda_{n_k}\in E_{n_k}$, for all $k\in\mathbb{N}$, such that $\lim \lambda_{n_k}=\lambda$.
\item $\lambda\in \liminf E_n$ if and only if there exists a sequence  $\{\lambda_n\}$  such that $\lambda_n\in E_n$ for all $n\in\mathbb{N}$, and $\lim \lambda_n = \lambda$.
\end{enumerate}
\end{rem}

\begin{prop}\label{prop0}
Let  $K,E,E_n$ be non-empty compact subsets of $\mathbb{C}$  such that $E_n\subseteq K$, for all $n\in\mathbb{N}$. Then $E_n\to E$ in the Hausdorff metric if and only if $\limsup E_n\subseteq E$ and
$E\subseteq\liminf E_n$.
\end{prop}

\section{Spectral continuity on complex Banach algebras}

Let $\mathcal{A}$ be a complex Banach algebra with identity $1_\mathcal{A}$. For $x\in\mathcal{A}$ the resolvent of $x$ is defined by $\rho(x)=\{\lambda\in\mathbb{C}: x-\lambda1_\mathcal{A} \textrm{ is invertible in }\mathcal{A}\}$ and the spectrum of $x$ is given by $\sigma(x)=\mathbb{C}\setminus\rho(x)$.  The spectral radius $r(x)$ of  $x$ is the number $r(x)=\max\{|\lambda|: \lambda\in\sigma(x)\}$. It is well known that $r(x)=\underset{n\to\infty}{\lim}\|x^n\|^\frac{1}{n}=\underset{n}{\inf}\|x^n\|^\frac{1}{n}$.

 A sequence $(x_n)$ in $\mathcal{A}$ is said to be norm convergent to $x$ (in notation $x_n\overset{n}{\to} x$), if $\|x_n-x\|\to 0$. M.  Ahues in \cite{MR1886113} introduces a new mode of convergence on $B(X)$, named $\nu-$convergence. This type of convergence can be generalized in the same way to complex unital Banach algebras: a sequence $(x_n)$ in $\mathcal{A}$ is said to be $\nu-$convergent to $x$, denoted by $x_n\overset{\nu}{\to}x$, if $(\|x_n\|)$ is bounded, $\|(x_n-x)x\|\to 0$ and $\|(x_n-x)x_n\|\to0$. The $\nu$-convergence is a pseudo-convergence in the sense that it is possible to have $x_n\overset{\nu}{\to}x$ and $x_n\overset{\nu}{\to}y$ but $x\neq y$, see for instance \cite[Example 1]{MR1886113}. There is a connection between norm convergence and $\nu$-convergence as follows: if $x_n\overset{n}{\to} x$ then $x_n\overset{\nu}{\to}x$, also, if $x_n\overset{\nu}{\to}x$ and $x$ is right invertible then $x_n\overset{n}{\to} x$. 
 
A function $\tau$ defined on $\mathcal{A}$ whose values are non-empty compact sets in $\mathbb{C}$ is said to be continuous ($\nu$-continuous) at $x$, if $\tau(x_n)\to\tau(x)$ with respect to the Hausdorff metric, for all sequence $(x_n)$ in $\mathcal{A}$ such that $x_n\overset{n}{\to}x$ ($x_n\overset{\nu}{\to}x$). It is clear that if $\tau$ is $\nu$-continuous at $x$, then $\tau$ is continuous at $x$.   $\tau$ is said to be upper semi-continuous (upper semi-$\nu$-continuous) at $x$, if $\limsup\tau(x_n)\subseteq \tau(x)$  for all sequence $(x_n)$ in $\mathcal{A}$ such that $x_n\overset{n}{\to}x$ ($x_n\overset{\nu}{\to}x$). Also, $\tau$  is said to be lower semi-continuous (lower semi-$\nu$-continuous) at $x$, if $\tau(x)\subseteq\liminf\tau(x_n)$  for all  $(x_n)$ in $\mathcal{A}$ such that $x_n\overset{n}{\to}x$ ($x_n\overset{\nu}{\to}x$). From Proposition \ref{prop0} we have next remark.

\begin{rem}\label{rem0}
Let $\tau$ be a function defined on $\mathcal{A}$ whose values are non-empty compact sets in $\mathbb{C}$ such that $\tau(y)\subseteq\sigma(y)$ for all $y\in\mathcal{A}$. Then 
\begin{enumerate}
\item $\tau$ is continuous  at $x\in\mathcal{A}$ if and only if $\tau$ is upper and lower semi-continuous at $x$.
\item $\tau$ is $\nu$-continuous  at $x\in\mathcal{A}$ if and only if $\tau$ is upper and lower semi-$\nu$-continuous at $x$.
\end{enumerate}
\end{rem}

Next theorem is proved in \cite[Corollary 2.7]{MR1886113}.

\begin{thm}\label{prop2a}
For each  $x\in\mathcal{A}$, $\sigma$ is upper semi-$\nu$-continuous at $x$.
\end{thm}

In the following proposition we will use the notation $\sigma_\mathcal{B}(x)$ for  the spectrum of $x\in\mathcal{B}$ with respect to  a given subalgebra $\mathcal{B}$ of $\mathcal{A}$. A character on an  abelian Banach algebra $\mathcal{B}$ is a non-zero homomorphism $\varphi:\mathcal{B}\to\mathbb{C}$. The set of all characters on $\mathcal{B}$ is denoted by $\mathcal{M}(\mathcal{B})$. It is well known that if $\mathcal{B}$ is an unital abelian Banach algebra then
\begin{enumerate}
\item $\varphi(x)\in\sigma_\mathcal{B}(x)$ for all  $x\in\mathcal{B}$ and $\varphi\in\mathcal{M}(\mathcal{B})$;
\item for every $\lambda\in\sigma_\mathcal{B}(x)$, there exists  $\varphi\in \mathcal{M}(\mathcal{B})$ such that $\varphi(x)=\lambda$;
\item for each $\varphi\in\mathcal{M}(\mathcal{B})$, $\|\varphi\|=1$. 
\end{enumerate}

A very well known result in spectral continuity is that if the elements of a sequence $(x_n)$ in a Banach algebra $\mathcal{A}$ commute with their limit, i.e. if $x_n\overset{n}{\to}x$ and $x_nx=xx_n$, then $\sigma(x_n)\to\sigma(x)$. We generalize this result  for the $\nu$-convergence, as shown in the following proposition.

\begin{prop}
Let $a\in\mathcal{A}$ and $(a_n)$ be a sequence in $\mathcal{A}$ such that $a_na=aa_n$ and $a_na_m=a_ma_n$ for all $n,m\in\mathbb{N}$. If $a_n\overset{\nu}{\to}a$ and $0\in\textup{acc\,}\sigma(a)$, then $\sigma(a_n)\to\sigma(a)$.
\end{prop}
\begin{proof}
By Theorem \ref{prop2a}, $\sigma$ is upper semi-$\nu$-continuous at $a$, thus $\limsup\sigma_\mathcal{A}(a_n)\subseteq\sigma_\mathcal{A}(a)$. Hence by Proposition \ref{prop0}, we only need to prove that $\sigma_\mathcal{A}(a)\subseteq\liminf\sigma_\mathcal{A}(a_n)$. 

Consider $\mathcal{B}_0$ the subalgebra of $\mathcal{A}$ which consists of all  linear combinations of finite products of elements in $\{a_n:n\in\mathbb{N}\}\cup\{a,1_\mathcal{A}\}$. From hypothesis, $\mathcal{B}_0$ is commutative. Thus by Zorn's lemma, there exists $\mathcal{B}$  the maximal abelian  subalgebra of $\mathcal{A}$ such that $\mathcal{B}_0\subseteq\mathcal{B}$.  Therefore by \cite[Exercise 8, p. 8]{MR1074574},  $\sigma_\mathcal{A}(a)=\sigma_\mathcal{B}(a)$ and $\sigma_\mathcal{A}(a_n)=\sigma_\mathcal{B}(a_n)$ for all $n\in\mathbb{N}$. Let $\lambda\in\sigma_\mathcal{A}(a)$ with $\lambda\neq 0$. Then there exists $\varphi\in \mathcal{M}(\mathcal{B})$ such that $\varphi(a)=\lambda$. Observe that $$\|(\varphi(a_n)-\varphi(a))\varphi(a)\|=\|\varphi((a_n-a)a)\|\leq\|\varphi\|\|(a_n-a)a\|\to0.$$
Thus $\|(\varphi(a_n)-\varphi(a))\lambda\|\to0$ which implies that $\varphi(a_n)\to\lambda$. Now, since $\varphi(a_n)\in\sigma_\mathcal{B}(a_n)(=\sigma_\mathcal{A}(a_n))$ for all $n\in\mathbb{N}$, it follows from Remark \ref{rem2} that $\lambda\in\liminf\sigma_\mathcal{A}(a_n)$. Consequently, since $0\in\textup{acc\,}\sigma_\mathcal{A}(a)$, $$\sigma_\mathcal{A}(a)=\overline{\sigma_\mathcal{A}(a)\setminus\{0\}}\subseteq\liminf\sigma_\mathcal{A}(a_n).$$
\end{proof}

An elementary Cauchy domain is a bounded open connected subset of $\mathbb{C}$ whose boundary is the union of a finite number of nonintersecting Jordan curves. A finite union of elementary Cauchy domains having disjoint closures is called a Cauchy domain. Let $D$ be a Cauchy domain, if each curve involved in the boundary of $D$ is oriented in such a way that points in $D$ lie to the left as the curve is traced out, then the oriented boundary $C$ of $D$ is called a Cauchy contour. The interior of the Cauchy contour $C$ is defined as $\textrm{int}(C)=D$ and the exterior of $C$ is defined by $\textrm{ext}(C)=\mathbb{C}\setminus(D\cup C)$.

 A set $\Lambda\subseteq\sigma(a)$ is a spectral set for $a$ if $\Lambda$ is closed as well as open in $\sigma(a)$.  We set $\mathcal{C}(a,\Lambda)$ the set of all Cauchy contours $C$ separating $\Lambda$ from $\sigma(a)\setminus\Lambda$, i.e. $\Lambda\subseteq \textrm{int}(C)$ and $\sigma(a)\setminus\Lambda\subseteq \textrm{ext}(C)$. It is clear that if $C\in \mathcal{C}(a,\Lambda)$ then $C\subseteq \rho(a)$. For a spectral set $\Lambda$ for $a$ and $C\in \mathcal{C}(a,\Lambda)$, define $$p(a,\Lambda)=-\frac{1}{2\pi i}\int_C(a-z)^{-1}dz.$$ The element $p(a,\Lambda)$ does not depend on the choice of $C\in \mathcal{C}(a,\Lambda)$.
\begin{rem}\label{remc}
Let $a\in \mathcal{A}$, $\Lambda$ a spectral set for $a$ and $C\in \mathcal{C}(a,\Lambda)$. If $p=p(a,\Lambda)$, then
\begin{enumerate}
\item $p^2=p$ and $pa=ap$
\item $\Lambda=\emptyset$  if and only if $p(a,\Lambda)=0$.
\end{enumerate}
\end{rem}
 \begin{prop}\label{prop1}
Let $p,q\in\mathcal{A}$ be such that $p^2=p$. If $p\neq 0$ and $r(p-q)<1$ then $q\neq 0$.
 \end{prop}
 \begin{proof}
 Since $r(p-q)<1$ it follows that $(p-q)-1_\mathcal{A}$ is invertible. Suppose that $q=0$, then $p-1_\mathcal{A}$ is invertible. Thus there exists $z\in\mathcal{A}$ such that $(p-1_\mathcal{A})z=1_\mathcal{A}$. This implies that $p(p-1_\mathcal{A})z=p1_\mathcal{A}$ and so $(p^2-p)z=p$ i.e. $0=(p-p)z=p$, which is a contradiction.
 \end{proof}
\begin{thm}\cite[Proposition 2.9]{MR1886113}\label{thm1}
Let $a\in\mathcal{A}$, $\Lambda$ be a spectral set for $a$,  $C\in\mathcal{C}(a,\Lambda)$ and $(a_n)$ be a sequence in $\mathcal{A}$ such that $a_n\overset{\nu}{\to}a$. Then
\begin{enumerate}
\item There exists $n_0\in\mathbb{N}$ such that for every $n\geq n_0$, $C$ lies in $\rho(a_n)$. 
\item If  $\Lambda_n:=\sigma(a_n)\cap \textup{int}(C)$  for all $n\geq n_0$, then $\Lambda_n$ is a spectral set for $a_n$ and $C\in \mathcal{C}(a_n,\Lambda_n)$. Further, if $0\in \textup{ext}(C)$ then  $p(a_n,\Lambda_n)\overset{\nu}{\to}p(a,\Lambda)$.
\end{enumerate}
\end{thm}

Next lemma is a generalization of \cite[Lemma 1.5]{MR543882}.
\begin{lem}\label{lemmma1}
Let $a\in \mathcal{A}$ and $(a_n)$ be a sequence in  $\mathcal{A}$ such that  $a_n\overset{\nu}{\to} a$. If $U$ is an open set  for which $0\not\in U$ and $U$ contains a component of $\sigma(a)$, then there exits  $n_0\in\mathbb{N}$ such that  $U$ contains a component of  $\sigma(a_n)$ for all $n\geq n_0$.
\end{lem}
\begin{proof}
Let $\Omega$ be a component of $\sigma(a)$ and  $U$ be an open set of $\mathbb{C}$ such that $0\not\in U$ and $\Omega\subseteq U$. Since  $\Omega\cap [\sigma(a)\setminus U]=\emptyset$, $\sigma(a)\setminus U$ is closed and $\sigma(a)$ is compact, it follows that there exists  $\Lambda\neq\emptyset$ open and closed set in $\sigma(a)$ such that $\Omega\subseteq \Lambda$ and $\Lambda\cap [\sigma(a)\setminus U]=\emptyset$. This implies that $\Lambda\subseteq U$.  From \cite[Theorem 1.21]{MR1886113}, there exists a Cauchy domain $D$ such that 
\begin{equation}\label{equation1}
\Lambda\subset D\hspace{.2cm}\textrm{ and }\hspace{.2cm} \overline{D}\subset U\cap [\mathbb{C}\setminus(\sigma(a)\setminus \Lambda)].
\end{equation}

Let $C$ be the boundary of $D$ oriented in way that $C$ be a Cauchy contour. It is clear by (\ref{equation1}) that $C\in\mathcal{C}(a,\Lambda)$. Then from Theorem \ref{thm1}, there exits $n_1\in\Bbb{N}$ such that for every $n\geq n_1$, 	$C$ lies in $\rho(a_n)$. Further, if $\Lambda_n:=\sigma(a_n)\cap D$ for all $n\geq n_1$ then since $0\not\in U$ we have that 
\begin{equation*}
[r(p-p_n)]^2\leq \|(p-p_n)^2\|\leq \|(p_n-p)p_n\|+\|(p_n-p)p\|\to 0,
\end{equation*}
where $p_n=p(a_n,\Lambda_n)$ and $p=p(a,\Lambda)$. Thus there exists $n_0\in\mathbb{N}$ with $n_0\geq n_1$ such that $r(p-p_n)<1$ for all $n\geq n_0$.  Since $\Lambda\neq\emptyset$ we have from Remark \ref{remc} that $p\neq 0$. Therefore, by Proposition \ref{prop1}, $p_n\neq 0$ for all $n\geq n_0$. Thus by Remark \ref{remc}, $\Lambda_n\neq \emptyset$ for all $n\geq n_0$.  This implies, since $\Lambda_n$ is both open and closed  in $\sigma(a_n)$, that there exists $\Omega_n$ a component of $\sigma(a_n)$ such that $\Omega_n\subseteq\Lambda_n$. Observe that $\Lambda_n\subseteq D \subseteq U$. Thus $\Omega_n\subseteq U$. Therefore $U$ contains a component of  $\sigma(a_n)$ for all $n\geq n_0$.
\end{proof}

\begin{rem}
From Lemma \ref{lemmma1} we have that if $\lambda\in\textup{iso\,}\sigma(a)$ with $\lambda\neq 0$, then $\lambda\in\liminf\sigma(a_n)$ for all sequence $(a_n)$ in $\mathcal{A}$ such that $a_n\overset{\nu}\to a$.
\end{rem}

\section{Spectral continuity on the algebra $B(X)$}
Let $X$ be a Banach space and let $B(X)$ be the algebra of all bounded linear operators defined on $X$. For $T\in B(X)$, let $N(T)$ and $R(T)$ be denote the null space and the range of the mapping $T$. Let $\alpha(T) = \dim N(T)$ and $\beta(T) = \dim X/R(T)$, if these spaces are finite dimensional, otherwise let $\alpha(T)=\infty$ and $\beta(T)=\infty$. If the range $R(T)$ of $T\in B(X)$ is closed and $\alpha(T)<\infty$ then $T$ is said to be an {\it upper semi-Fredholm} operator ($T\in \Phi_+(X)$). Similarly, if  $\beta(T)<\infty$ then $T$ is said to be a {\it lower semi-Fredholm} operator ($T\in\Phi_-(X)$). If $T\in \Phi_-(X)\cup\Phi_+(X)$ then $T$ is called a
{\it semi-Fredholm} operator ($T\in \Phi_{\pm}(X)$) and for  $T\in \Phi_-(X)\cap\Phi_+(X)$ we say that  $T$ is  a {\it Fredholm} operator ($T\in \Phi(X)$). For  $T\in \Phi_{\pm}(X)$, the {\it index} of $T$ is defined by $i(T)=\alpha(T)-\beta(T).$ It is well known that the index is a continuous function on the set of semi-Fredholm operators. This property also holds for the $\nu$-convergence in the following sense:

\begin{thm}\cite[Theorem 3.4]{MR3388799}\label{thmindex}
Let $T\in \Phi(X)$ and $(T_n)$ be a sequence in $B(X)$ such that $(T_n-T)T\overset{n}{\to}0$. If $\lambda,\lambda_n$ are complex numbers such that $\lambda_n\to\lambda$ and $$T-\lambda,T_n-\lambda_n\in \Phi_\pm(X)$$ for all $n\in\mathbb{N}$, then $i(T_n-\lambda_n)\to i(T-\lambda)$.
\end{thm}

 For $T\in B(X)$ and $k\in\mathbb{N}\cup\{-\infty,\infty\}$, let $\rho_{sf}^k(T)$ be  denote  the set of $\lambda\in\mathbb{C}$ for which  $ T-\lambda\in\Phi_{\pm}(X)$ and  $i(T-\lambda )=k$. Put 
\begin{align*}
\rho^{-}_{sf}(T)=\underset{-\infty\leq k\leq -1}{\cup}\rho_{sf}^k(T),\hspace{.3cm}\rho^{+}_{sf}(T)=\underset{1\leq k\leq\infty}{\cup}\rho_{sf}^k(T),\hspace{.3cm}\rho_{sf}^\pm(T)=\rho_{sf}^-(T)\cup\rho_{sf}^+(T).
\end{align*}
All these sets are open and are contained in $\sigma(T)$. Let $\sigma_{sf}(T)$ be denote the set of all $\lambda\in \mathbb{C}$ such that $T-\lambda\not\in\Phi_\pm(X)$.
Let $K(X)$ be denote the set of all compact linear operators in $B(X$). If $\pi : B(X) \to B(X)/K(X)$ is the canonical homomorphism, then the essential spectrum of an operator $T \in B(X)$, $\sigma_e(T)$, is the spectrum of $\pi(T)$ in the Calkin algebra $B(X)/K(X)$. Also, the left essential spectrum $\sigma_{le}(T)$ (the right essential spectrum $\sigma_{re}(T)$) is the left spectrum (right spectrum) of $\pi(T)$. We set $\sigma_{lre}(T)=\sigma_{le}(T)\cap\sigma_{re}(T)$.  It is clear that 
\begin{equation}\label{eq0}
\sigma_{sf}(T)\subseteq\sigma_{lre}(T)
\end{equation}
but the opposite inclusion is not always satisfied in general Banach spaces. These classes of sets coincide in the case of Hilbert spaces.

The approximate point spectrum, the surjective spectrum, the point spectrum, the Weyl spectrum  and the set Riesz points of $T\in B(X)$ are defined by $\sigma_{ap}(T)=\{\lambda\in\mathbb{C}:T-\lambda \textrm{ is not bounded below}\}$, $\sigma_s(T)=\{\lambda\in\mathbb{C}:T-\lambda \textrm{ is not surjective}\}$, $\sigma_p(T)=\{\lambda\in\mathbb{C}:\lambda$ is an eigenvalue of $T\}$, $\sigma_w(T)=\{\lambda\in\mathbb{C}:T-\lambda$ is not a Fredholm operator of index zero$\}$ and $\pi_0(T)=\{\lambda\in\mathbb{C}: \lambda$  is an isolated eigenvalue of $T$ of finite algebraic multiplicity$\}$.

\begin{rem}\label{rem1}
Let $T\in B(X)$  and $(T_n)$ be a sequence in $B(X)$ such that $T_n\overset{\nu}{\to} T$. The following inclusions are hold:
\begin{enumerate}
\item $\pi_0(T)\subseteq\liminf\pi_0(T_n)$. See, \cite[Corollary 2.13]{MR1886113}.
\item $\Delta\subseteq\liminf\sigma(T_n)$, where $\Delta=\big\{\lambda\in\sigma(T)\setminus(\overline{\rho_{sf}^\pm(T)}\cup\{0\}): $ for every $\epsilon>0$, the ball $B(\lambda,\epsilon)$ contains a component of $\sigma_{sf}(T)\cup\pi_0(T)\big\}$ and $T$ satisfies Browder's theorem. Really, take $\lambda\in \Delta$ then $\lambda\in\sigma(T)\setminus\overline{\rho_{sf}^\pm(T)}$ and $\lambda\neq 0$, so there exists $r>0$ such that $B(\lambda,r)\cap \overline{\rho_{sf}^\pm(T)}=\emptyset$ and $0\not\in B(\lambda,r)$. Let $\epsilon>0$ with $\epsilon<r$, since $\lambda\in\Delta$, the ball $B(\lambda,\epsilon)$ contains a component $\Omega$ of $\sigma_{sf}(T)\cup\pi_0(T)$. By \cite[Lemma 3.6]{sanchez4}, $\Omega$ is a component of $\sigma(T)$. Therefore by Lemma \ref{lemmma1}, there exists $n_0\in\mathbb{N}$ such that $B(\lambda,\epsilon)$ contains a component $\Omega_n$ of $\sigma(T_n)$ for all $n\geq n_0$. Thus $B(\lambda,\epsilon)\cap \sigma(T_n)\neq \emptyset$ for all $n\geq n_0$. This implies that $\lambda\in \liminf\sigma(T_n)$.
\item $[\textup{iso}\,\sigma(T)]\setminus\{0\}\subseteq\liminf\sigma(T_n)$.
\suspend{enumerate}
Moreover, if $T\in \Phi(X)$ then
\resume{enumerate}
\item $\overline{\rho_{sf}^+(T)}\subseteq\liminf\sigma_a(T_n)$. See, \cite[Theorem 3.6]{MR3388799}.
\item $\overline{\rho_{sf}^-(T)}\subseteq\liminf\sigma_s(T_n)$.
\item $\overline{\rho_{sf}^\pm(T)}\subseteq\liminf\sigma_w(T_n)$.
\item $\Gamma\subseteq\liminf\sigma_{sf}(T_n)$, where $\Gamma=\big\{\lambda\in\sigma(T):$ for every $\epsilon>0$, there exist points $\mu_1,\mu_2\in B(\lambda,\epsilon)$ such that $T-\mu_1,T-\mu_2\in \Phi_\pm(X)$ and $i(T-\mu_1)\neq i(T-\mu_2)\big\}$. This result was established in \cite[Proposition 2.1]{MR634170} for the norm convergence. The proof for the $\nu$-convergence is in a similar way.
\end{enumerate}
\end{rem}

\begin{thm}\label{thm3.2}
Let $T\in \Phi(X)$ be such that $0\in\textup{acc}\sigma(T)$. If for each $\lambda\in\sigma(T)\setminus\overline{\rho_{sf}^\pm(T)}$ and $\epsilon>0$, the ball $B(\lambda,\epsilon)$ contains a component of $\sigma_{sf}(T)\cup\pi_0(T)$, then $\sigma$ is $\nu$-continuous at $T$.
\end{thm}
\begin{proof}
First observe that $\sigma(T)\setminus\sigma_w(T)=\pi_0(T)\cup\textrm{int}[\sigma(T)\setminus\sigma_w(T)]$, see for example \cite[Proposition 2.1]{MR2894894}. By hypothesis we have that $\textrm{int}[\sigma(T)\setminus\sigma_w(T)]=\emptyset$, thus $T$ satisfies Browder's theorem. 

From Remark \ref{rem0} and Theorem \ref{prop2a} it is sufficient to prove that $\sigma$ is lower semi-$\nu$-continuous at $T$.  Let $(T_n)$ be a sequence in $B(X)$ such that $T_n\overset{\nu}{\to}T$.  By Remark \ref{rem1} (6), $\overline{\rho_{sf}^\pm(T)}\subseteq\liminf\sigma(T_n)$. Now, from hypothesis, $\sigma(T)\setminus(\overline{\rho_{sf}^\pm(T)}\cup\{0\})\subseteq\Delta$, therefore by Remark \ref{rem1} (2), $\sigma(T)\setminus(\overline{\rho_{sf}^\pm(T)}\cup\{0\})\subseteq\liminf\sigma(T_n)$.  Consequently, $\sigma(T)=\overline{\sigma(T)\setminus\{0\}}\subseteq\liminf\sigma(T_n)$. Thus $\sigma$ is lower semi-$\nu$-continuous at $T$. 
\end{proof}

\begin{cor}
Let $H$ be a Hilbert space and $T\in \Phi(H)$ be such that $0\in\textup{acc}\sigma(T)$. Then, $\sigma$ is continuous at $T$ if and only if $\sigma$ is $\nu$-continuous at $T$.
\end{cor}
\begin{proof}
It is clear that the $\nu$-continuity of $\sigma$ at $T$ implies the continuity of $\sigma$ at $T$. Now, if $\sigma$ is continuous at $T$, then by \cite[Theorem 3.1]{MR543882}, for each $\lambda\in\sigma(T)\setminus\overline{\rho_{sf}^\pm(T)}$ and $\epsilon>0$, the ball $B(\lambda,\epsilon)$ contains a component of $\sigma_{sf}(T)\cup \pi_0(T)$. Therefore, by Theorem \ref{thm3.2}, $\sigma$ is $\nu$-continuous at $T$.
\end{proof}
\begin{thm}\label{limsupap}
Let  $T\in \Phi(X)$. If  $(T_n)$ is a sequence in $B(X)$ such that $T_n\overset{\nu}{\to}T$ then $$[\lim\sup \sigma_{ap}(T_n)]\setminus\{0\}\subseteq\sigma_{ap}(T).$$
\end{thm}
\begin{proof}
Let $D,E$ be closed subspaces of $X$ with $\dim E<\infty$ such that 
\begin{equation}\label{sum}
X=N(T)\oplus D \hspace{.2cm}\textrm{ and }\hspace{.2cm} X=R(T)\oplus E.
\end{equation}

Consider $(T_n)$  a sequence in $B(X)$ which is $\nu-$convergent to $T$. Let $\lambda\in \limsup\sigma_{ap}(T_n)$ with $\lambda\neq 0$. By Remark \ref{rem2}, there exist an increasing sequence of natural numbers $(n_k)$ and points $\lambda_{n_k}\in\sigma_{ap}(T_{n_k})$ such that $\lambda_{n_k}\to\lambda$. Suppose that $\lambda\not \in \sigma_{ap}(T)$. Then $T-\lambda\in\Phi_+(X)$ and $\alpha(T-\lambda)=0$. By (\ref{sum}), $R(T|_D)=R(T)$ and so $T|_D$ is bounded below, therefore by \cite[Theorem 5.26]{MR1861991}, $(T-\lambda)T|_D\in \Phi_+(D,X)$ and $\alpha((T-\lambda)T|_D)=0$.

On the other hand, observe that 
\begin{align*}
(T_{n_k}-\lambda_{n_k})T|_D&=(T-\lambda)T|_D+(T_{n_k}-T)T|_D+(\lambda-\lambda_{n_k})T|_D.
\end{align*}
From $\|(T_n-T)T\|\to 0$ we have that  $(T_{n_k}-\lambda_{n_k})T|_D$ converges in norm to $(T-\lambda)T|_D$. Consequently  by \cite[Theorem 5.23]{MR1861991}, there exists $k_0\in\mathbb{N}$ such that every $k\geq k_0$, $(T_{n_k}-\lambda_{n_k})T|_D\in \Phi_+(D,X)$ and $\alpha((T_{n_k}-\lambda_{n_k})T|_D)=0$. 

Suppose that for each $k\geq k_0$, $$N(T_{n_k}-\lambda_{n_k})\cap E\neq\{0\}.$$ Take $v_k\in N(T_{n_k}-\lambda_{n_k})\cap E$ with $\|v_k\|= 1$ for all $k\geq k_0$. Since $\dim E<\infty$ it follows that $F:=\{e\in E:\|e\|=1\}$ is a compact set. Therefore we may assume without loss of generality that there exists $v\in F$ such that $v_{k}\to v$. Observe that
\begin{align*}
\|(T_{n_k}-T)T_{n_k}\|&\geq \|(T_{n_k}-T)T_{n_k}v_{k}\|\\
&=\|(T_{n_k}-T)\lambda_{n_k}v_{k}\|\\
&=|\lambda_{n_k}|\|\lambda_{n_k}v_k-Tv_k\|
\end{align*}
for all $k\geq k_0$, and $ |\lambda_{n_k}|\|\lambda_{n_k}v_k-Tv_k\|\to |\lambda|\|\lambda v-Tv\|$.  This implies that $$|\lambda|\|\lambda v-Tv\| = \lim|\lambda_{n_k}|\|\lambda_{n_k}v_k-Tv_k\|\leq \lim \|(T_{n_k}-T)T_{n_k}\|=0,$$ and so $\|\lambda v-Tv\|=0$, i.e. $Tv=\lambda v$. Consequently, $\lambda\in\sigma_p(T)(\subseteq\sigma_{ap}(T))$, which is a contradiction. Therefore there exists $k'\geq k_0$ such that $N(T_{n_{k'}}-\lambda_{n_{k'}})\cap E=\{0\}$. Thus
\begin{equation}\label{suma2}
X=R(T)\oplus N(T_{n_{k'}}-\lambda_{n_{k'}})\oplus E.
\end{equation}
Then by (\ref{sum}) and (\ref{suma2}), 
\begin{align*}
\dim E=\dim X/R(T)& =\dim[N(T_{n_{k'}}-\lambda_{n_{k'}})\oplus E]\\
&=\dim N(T_{n_{k'}}-\lambda_{n_{k'}})+\dim E.
\end{align*}
Hence $\dim  N(T_{n_{k'}}-\lambda_{n_{k'}})=0$ and so $N(T_{n_{k'}}-\lambda_{n_{k'}})=\{0\}$. From (\ref{suma2}), $$R(T_{n_{k'}}-\lambda_{n_{k'}})=(T_{n_{k'}}-\lambda_{n_{k'}})T(D)+(T_{n_{k'}}-\lambda_{n_{k'}})(E),$$ which implies that $R(T_{n_{k'}}-\lambda_{n_{k'}})$ is closed. Therefore $\lambda_{n_{k'}}\not\in\sigma_{ap}(T_{n_{k'}})$, a contradiction. Consequently, $\lambda\in \sigma_{ap}(T)$.
\end{proof}

Let $\phi_+(T)$ be the set of $\lambda\in\rho^{+}_{sf}(T)$  such that $N(T-\lambda)$  is complemented in $X$ and $\phi_-(T)$ the set of $\lambda\in \rho^-_{sf}(T)$ such that $R(T-\lambda)$ is complemented in $X$. We set $\phi_\pm(T)=\phi_+(T)\cup\phi_-(T)$,  $\phi_{+\infty}(T)=\{\lambda\in\phi_+(T)\ |\ i(T-\lambda)=\infty\}$,  $\phi_{-\infty}(T)=\{\lambda\in\phi_-(T)\ |\ i(T-\lambda)=-\infty\}$ and $\phi_{\pm\infty}(T)=\phi_{+\infty}(T)\cup\phi_{-\infty}(T)$.  Al these sets are open and it is clear that 
\begin{equation}
\sigma_e(T)=\sigma_{lre}(T)\cup\phi_{\pm\infty}(T).
\end{equation}

With this equality, \cite[Proposition 1.3]{MR543882} can be extended to general Banach spaces:
\begin{prop}\label{prop2}
If $C$ is a component of $\sigma_{lre}(T)$ and $C\cap\overline{\phi_\pm(T)}=\emptyset$ then $C$ is a component of $\sigma_e(T)$.
\end{prop}
\begin{thm}
If the following conditions hold
\begin{enumerate}[i.]
\item $0\in\textup{acc}\sigma_{ap}(T)$
\item $T\in \Phi(X)$ 
\item $T$ satisfies Browdwer's theorem
\item $\sigma_p(T)\cap\phi_{-}(T)=\emptyset$
\item for every $\epsilon>0$ and $\lambda\in \sigma_{lre}(T)\setminus\overline{\phi_\pm(T)}$, the ball $B(\lambda,\epsilon)$ contains a component of $\sigma_{lre}(T)$,
\end{enumerate}
then $\sigma_{ap}(T_n)\to\sigma_{ap}(T)$ for all $T_n\overset{\nu}{\to}T$.
\end{thm}
\begin{proof}
From conditions (i) and (ii) we have by Theorem \ref{limsupap} that $$\limsup\sigma_{ap}(T_n)\subseteq\sigma_{ap}(T).$$

Now, we show that $\sigma_{ap}(T)\subseteq\liminf\sigma_{ap}(T_n)$. Take $\lambda\in\sigma_{ap}(T)$ with $\lambda\neq 0$. Suppose that $\lambda\in\sigma_{lre}(T)\setminus\overline{\phi_\pm(T)}$. Let $\epsilon>0$, there exists $0<\epsilon_1<\epsilon$ such that $B(\lambda,\epsilon_1)\cap\overline{\phi_\pm(T)}=\emptyset$ and $0\not\in B(\lambda,\epsilon_1)$. By  hypothesis (v), $B(\lambda,\epsilon_1)$ contains a component  $\Omega$  of $\sigma_{lre}(T)$. Then  by Proposition \ref{prop2}, $\Omega$ is a component of $\sigma_e(T)$. Since $T_n\overset{\nu}{\to}T$ we have that $\pi(T_n)\overset{\nu}{\to}\pi(T)$ in the Calkin algebra $B(X)/K(X)$. Thus by Lemma \ref{lemmma1}, there exists $n_0\in\mathbb{N}$ such that $B(\lambda,\epsilon_1)$ contains a component $\Omega_n$ of $\sigma_e(T_n)$ for all $n\geq n_0$. Consequently $$\emptyset\neq \partial \Omega_n\subseteq\partial \sigma_e(T_n)\cap B(\lambda,\epsilon)\subseteq \sigma_{sf}(T_n)\cap B(\lambda,\epsilon)$$
for all $n\geq n_0$. Therefore $\lambda\in\liminf\sigma_{sf}(T_n)(\subseteq\liminf\sigma_{ap}(T_n))$.

If $\lambda\not\in\sigma_{lre}(T)\setminus\overline{\phi_\pm(T)}$, then $\lambda\in\phi_\pm(T)$ or $i(T-\lambda)=0$. From condition (iv), we have that $\lambda\not\in\phi_-(T)$. Thus from condition (iii), $\lambda\in  \phi_+(T)\cup\pi_0(T)$. Therefore, by Remark \ref{rem1}, $\lambda\in\liminf\sigma_{ap}(T_n)$.
\end{proof}

\section{On certain class of operators}

We say that an operator $T\in B(H)$ is essentially $G_1$ if $\|(\pi(T)-z)^{-1}\|=\frac{1}{d(z,\sigma_e(T))}$ for all $z\not\in \sigma_e(T)$. In \cite[Theorem 6]{MR433257} it is shown that the restriction of the Wey spectrum  on the class of essentially $G_1$ operators is continuous. This is also true for $\nu$-continuity, as the following theorem states.

\begin{thm}\label{thmB}
Let $T\in \Phi(H)$ be such that $0\not\in \sigma_w(T)$ or $0\in\textrm{acc}\sigma_w(T)$. If $(T_n)$ is a sequence of  essentially $G_1$ operators such that $T_n\overset{\nu}{\to}T$ then $\sigma_w(T_n)\to\sigma_w(T)$.
\end{thm}
\begin{proof}
From  \cite[Theorem 3.3]{MR3624565} we have that $\limsup\sigma_w(T_n)\subseteq \sigma_w(T)$. We show that $\sigma_w(T)\subseteq\liminf\sigma_w(T_n)$. Let $\lambda\in\sigma_w(T)\setminus\{0\}$ and suppose that $\lambda\not\in\liminf\sigma_w(T_n)$.  Then, there exist $\epsilon>0$  such that  $B_\epsilon(\lambda)\cap\sigma_w(T_{n}) = \emptyset$ for  infinite number of $n$'s. Without loss of generality assume that this holds for all $n$. This implies that  $d(\lambda,\sigma_e(T_{n}))\geq d(\lambda,\sigma_w(T_{n}))\geq \epsilon>0$. Thus $\lambda\not\in\sigma_e(T_{n})$, now from the fact that $T_{n}$ is essentially $G_1$, we have that $\|(\pi(T_{n})-\lambda)^{-1}\|=\frac{1}{d(\lambda,\sigma_e(T_{n}))}\leq\frac{1}{\epsilon}$. This implies that

$\Big\|\Big[\big(\pi(T_n)-\pi(T)\big) (\pi(T_{n})-\lambda)^{-1}\Big]^2\Big\| \leq$
\begin{align}\label{eqvn0}
&\leq \frac{1}{|\lambda|}\Big[\|\big(\pi(T_n)-\pi(T)\big)\pi(T_{n})\|\|(\pi(T_{n})-\lambda)^{-1}\|\|\pi(T_n)-\pi(T)\|\nonumber\\ 
&\ \ \ \ \ +\ \|\big(\pi(T_n)-\pi(T)\big)\pi(T)\|+\|\big(\pi(T_n)-\pi(T)\big)\pi(T_n)\|\Big]\|(\pi(T_{n})-\lambda)^{-1}\|\nonumber\\
&\leq\frac{1}{|\lambda|}\Big[\|(T-T_{n})T_{n}\|\|T_n-T\|\frac{1}{\epsilon}+\|(T-T_{n})T\|+\|(T-T_{n})T_{n}\|\Big]\frac{1}{\epsilon}.
\end{align}

Now, since  $T_n\overset{\nu}{\to}T$, it follows that $\|(T_n-T)T\|\to 0$, $\|(T_n-T)T_{n}\|\to 0$ and $(\|T_{n}\|)$ is bounded. Therefore  the expression (\ref{eqvn0}) tends to 0. Thus, there exists $n^*\in\mathbb{N}$ such that $$\Big\|\Big[\big(\pi(T_{n^*})-\pi(T)\big) (\pi(T_{n^*})-\lambda)^{-1}\Big]^2\Big\|<1.$$

Then $r\Big(\big(\pi(T_{n^*})-\pi(T)\big) (\pi(T_{n^*})-\lambda)^{-1}\Big)<1,$ which implies that $\big(\pi(T_{n^*})-\pi(T)\big) (\pi(T_{n^*})-\lambda)^{-1}-1$ is invertible and so
\begin{align*}
\pi(T)-\lambda&=(\pi(T_{n^*})-\lambda)-(\pi(T_{n^*})-\pi(T))\\
&=-\Big[\big(\pi(T_{n^*})-\pi(T)\big) (\pi(T_{n^*})-\lambda)^{-1}-1\Big](\pi(T_{n^*})-\lambda)
\end{align*}
is invertible. Consequently, $\lambda\not\in\sigma_e(T)$, i.e. $T-\lambda$ is a Fredholm operator. Finally, since $\lambda\not\in\sigma_w(T_n)$ for all $n\in\mathbb{N}$, it follows by Theorem \ref{thmindex} that  $i(T_n-\lambda)\to i(T-\lambda)$, and so $i(T-\lambda)=0$, which is a contradiction.
\end{proof}

\begin{thm}
Let $T\in \Phi(H)$ be such that satisfies Browder's theorem. If $0\in\textrm{acc}\sigma_w(T)$ and $(T_n)$ is a sequence of  essentially $G_1$ operators such that $T_n\overset{\nu}{\to}T$, then $\sigma(T_n)\to\sigma(T)$.
\end{thm}

An operator $A\in B(H)$ is called $p$-hyponormal if $(A^*A)^p-(AA^*)^p\geq 0$. For the case $p=1$ the operator $A$ is called hyponormal. It is well known that the restriction of the spectrum on the class of $p$-hyponormal operators is continuous. See, \cite{MR1814478} and  \cite{MR1785076}. In the following theorem we use the idea of \cite{MR1814478} and adapt it for the case of $\nu$-convergence. Note first that if $A\in B(H)$ is a $p$-hyponormal operator and $0\in \sigma_p(A)$ then  $0\in \sigma_{jp}(A)$ and so $N(A)=N({A}^*)$ which implies that $N(A)$ is invariant for both $A$ and ${A}^*$. Therefore
\begin{equation}\label{eqbc}
A=\begin{bmatrix}
0&0\\
0&B
\end{bmatrix}
\end{equation}
on $N(A)\oplus N(A)^\perp$, where $B=A|_{N(A)^\perp}$ and $0\not\in\sigma_p(B)$. From \cite[Lemma 4]{MR1670405}, $B$ is also $p$-hyponormal. We claim that $0\not\in\sigma(|B|)$, indeed if $0\in\sigma(|B|)(=\sigma_{ap}(|B|))$ then there exists a sequence $(x_m)$ of unit vectors such that $|B|x_m\to0$. This implies that $Bx_m\to0$, thus $0\in\sigma_{ap}(B)$, but since $R(A)$ is closed i.e. $R(B)$ is closed, it follows that $0\in \sigma_p(B)$, which is a contradiction.

\begin{thm}\label{thmd}
If $T_n,T$ are operators in $B(H)$ such that 
\begin{enumerate}
\item $T_{n}\overset{\nu}{\to}T$,
\item $T\in \Phi(H)$ and $T_n$ is $p$-hyponormal for all $n\in\mathbb{N}$,
\item $0\in\sigma_p(T_n)$ for all $n\in\mathbb{N}$, and the sequence $(\||B_n|^{-1}\|)$ is bounded, where the operators $|B_n|$ are as in (\ref{eqbc}), 
\end{enumerate}
then $\sigma(T_{n})\to\sigma(T)$.
\end{thm}
\begin{proof}
First observe that $\|T_nT-T^2\|=\|(T_n-T)T\|\to0$, thus $T_nT\overset{n}{\to}T^2$, which implies that $T_nT\in \Phi(H)$ for all $n$ large. Thus we may suppose that $R(T_n)$ is closed for all $n\in\mathbb{N}$. We show that there exist a sequence $(S_n)$ of hyponormal operators and a sequence $(X_n)$ of invertible operators  such that $S_n=X_nT_nX^{-1}_n$ for all $n\in\mathbb{N}$, and $(\|X_{n}\|)$, $(\|X_{n}^{-1}\|)$ are  bounded. 

From condition (2), $0\in\sigma_{p}(T_{n})$ for all $n\in\mathbb{N}$. Then $T_{n}=\begin{bmatrix}
0&0\\
0&B_{n}
\end{bmatrix}$ on $N(T_{n})\oplus N(T_{n})^{\perp}$, $0\not\in\sigma_p(B_n)$, $B_n$ is $p$-hyponormal and $0\not\in\sigma(|B_n|)$. Consider the polar decomposition $B_n=U_n|B_n|$ and define $\widehat{B}_n=|B_n|^{1/2}U_n|B_n|^{1/2}$. Observe that if $x\in N(\widehat{B}_n)$ then $|B_n|^{1/2}U_n|B_n|^{1/2}x=0$ and so $B_n|B_n|^{-1}|B_n|^{1/2}x=0$ which implies that $|B_n|^{-1}|B_n|^{1/2}x=0$ because $0\not\in\sigma_p(B_m)$, hence $x=0$. Thus $N(\widehat{B}_n)=\{0\}$ i.e. $0\not\in\sigma_p(\widehat{B}_n)$. This implies that $0\not \in\sigma(|\widehat{B}_n|)$. Let $\widehat{B}_n$ have the polar decomposition $\widehat{B}_n=V_n|\widehat{B}_n|$, then by  \cite[Corollary 3]{MR1047771}, the operator $\widetilde{B}_n$, defined by $\widetilde{B}_n=|\widehat{B}_n|^{1/2}V_n|\widehat{B}_n|^{1/2}$, is hyponormal. Define
$$X_n=\begin{bmatrix}
1&0\\
0&|\widehat{B}_n|^{1/2}|B_n|^{1/2}
\end{bmatrix} \hspace{.3cm}\textrm{and}\hspace{.3cm} S_n=\begin{bmatrix}
0&0\\
0&\widetilde{B}_n
\end{bmatrix}.$$ Then  $S_n$ is hyponormal, $X_n$ is invertible and $S_n=X_nT_nX_n^{-1}$. From condition (2) we have that $(\|X_n^{-1}\|)$ is bounded. Also it is clear that $(\|X_n\|)$ is bounded.

We show that $\sigma(T)\subseteq\liminf\sigma(T_n)$. Take $\lambda\in\sigma(T)\setminus\{0\}$ and suppose that $\lambda\not\in\liminf\sigma(T_n)$. Then we may assume that there exists $\epsilon>0$ such that $B(\lambda,\epsilon)\cap \sigma(T_n)=\emptyset$ for all $n\in\mathbb{N}$. This implies that  $T_n-\lambda$ is invertible. In a similar way to proof of Theorem \ref{thmB}, we have that
\begin{align}\label{eqvn}
\Big\|\Big[\big(T_n-T) (T_{n}-\lambda)^{-1}\Big]^2\Big\| &\leq \frac{1}{|\lambda|}\Big[\|(T_n-T)T_{n}\|\|(T_{n}-\lambda)^{-1}\|\|T_n-T\|\nonumber\\ 
&\ \ \  +\|(T_n-T)T\|+\|(T_n-T)T_n\|\Big]\|(T_{n}-\lambda)^{-1}\|.
\end{align}

Since $T_n$ and $S_n$ are similar, it follows that $\sigma(T_n)=\sigma(S_n)$. Therefore, $d(\lambda,\sigma(S_n))$ $\geq \epsilon$ and $S_n-\lambda$ is invertible. Note that  $(T_n-\lambda)^{-1}=X_n^{-1}(S_n-\lambda)^{-1}X_n$. Moreover, since $S_n$ is hyponormal it follows that $\|(S_n-\lambda)^{-1}\|=\frac{1}{d(\lambda,\sigma(S_n))}\leq\frac{1}{\epsilon}$. Thus the right term of (\ref{eqvn})  is bounded by
\begin{equation}\label{eq13}
\frac{1}{|\lambda|}\Big[\|(T_n-T)T_{n}\|\frac{M_1M_2}{\epsilon}\|T_n-T\|+\|(T_n-T)T\|+\|(T_n-T)T_n\|\Big]\frac{M_1M_2}{\epsilon},
\end{equation}
where $M_1,M_2$ are constants such that $\|X_n\|\leq M_1$ and $\|X_n^{-1}|\leq M_2$ for all $n\in\mathbb{N}$. Since $T_n\overset{\nu}{\to}T$ it follows that the expresion in (\ref{eq13}) tends to zero. Proceeding similarly to the final part of the proof of Theorem \ref{thmB}, we obtain that $T-\lambda$ is invertible, which is a contradiction.
\end{proof}

\begin{rem}
The conclusion of Theorem \ref{thmd} holds if we replace the hypothesis by the following conditions:
\begin{enumerate}
\item $T_{n}\overset{\nu}{\to}T$, $T^*(T_n-T)\overset{n}{\to}0$ and $T_n^*(T_n-T)\overset{n}{\to}0$;
\item $T\in \Phi(H)$ and $\{0\}\neq N(T)\subseteq N(T_{n})$ for all $n\in\mathbb{N}$;
\item $T,T_n$ are $p$-hyponormal operators.
\end{enumerate}
\end{rem}
Indeed, from condition (1), $|T_n|^2-|T|^2=T_n^*T_n-T^*T=T_n^*(T_n-T)+[T^*(T_n-T)]^*\overset{n}{\to}0$,  thus $|T_n|^{1/2}\overset{n}{\to}|T|^{1/2}$. Since $0\in \sigma_p(T)$, it follows that $T=0\oplus B$ on $N(T)\oplus N(T)^{\perp}$ with $0\not\in\sigma(|B|)$. Then there exists $\alpha>0$ such that $\alpha\|y\|\leq \||B|^{1/2} y\|$ for all $y\in N(T)^\perp$.  This implies by condition (2) that for $0<\epsilon<\alpha$, there exists $N\in\mathbb{N}$ such that for every $n\geq N$, $(\alpha-\epsilon)\|y\|\leq \||B_n|^{1/2}y\|$ for all $y\in N(T_n)^\perp$. Consequently, $\|(|B_n|^{1/2})^{-1}\|\leq \frac{1}{\alpha-\epsilon}$ for all $n\geq N$.
\medskip

Berberian shows  that for every Hilbert space $H$, there exists a Hilbert space $K\supset H$ and a faithful $*-$representation $T\to T^\circ$ from $B(H)$ to $B(K)$: $(S+T)^\circ=S^\circ+T^\circ$, $(\lambda T)^\circ=\lambda T^\circ$, $(ST)^\circ = S^\circ T^\circ$, $(T^*)^\circ=(T^\circ)^*$, $(I_H)^\circ=I_K$ and $\|T^\circ\|=\|T\|$ such that
\begin{enumerate}
\item  $T\geq0$ if and only if $T^\circ\geq 0$,
\item $\sigma_p(T^\circ)=\sigma_a(T^\circ)=\sigma_a(T)$.
\end{enumerate}

\begin{rem} 
Observe that in the previous theorem, $\sigma_p(B_n)=\sigma_{ap}(B_n)$ due to $R(T_n)$ is closed. In \cite{MR1814478}  the authors  use the Berberian extension $T_n^\circ$ of a $p$-hyponormal operator  $T_n$ and state that if $0\in\sigma_p(T_n^\circ)$, then
\begin{equation}\label{equation14}
\sigma_p(B_n)=\sigma_{ap}(B_n),
\end{equation}
where $T_n^\circ=0\oplus B_n$ on $N(T_n)\oplus N(T_n)^\perp$ and $0\not\in \sigma_p(B_n)$, without the need for $R(T_n^\circ)$ to be closed. This fact was also established in \cite{MR2651679}, page 586, line 20. The authors claim that  
\begin{equation}\label{equation15}
\sigma_{ap}(B_\lambda)=\sigma_p(B_\lambda)
\end{equation}
for all  non-zero $\lambda\in\sigma_p(A^\circ),$ where  $A\in \mathcal{C}(i)$, and this collection is defined as the set of all operators $T\in B(H_i)$ for which $\sigma(T)=\{0\}$ implies $T$ is nilpotent and $T^\circ$ satisfies the property:
\begin{equation*}
 T^\circ =\begin{bmatrix}
\lambda&X_\lambda\\
0&B_\lambda
\end{bmatrix} \hspace{.2cm}\textrm{ on }\hspace{.2cm} N(T^\circ-\lambda)\oplus N(T^\circ-\lambda)^\perp
\end{equation*} 
at every non-zero $\lambda\in\sigma_p(T^\circ)$ for some operators $X_\lambda$ and $B_\lambda$ such that $\lambda\not\in\sigma_p(B_\lambda)$ and $\sigma(T^\circ)=\{\lambda\}\cup\sigma(B_\lambda)$.  However, equalities (\ref{equation14}) and (\ref{equation15}) are not necessary hold. We prove only that (\ref{equation15}) is false. It is clear that $\sigma_a(B_\lambda)\setminus\{\lambda\}=\sigma_p(B_\lambda)$ and $\alpha(B_\lambda-\lambda)=0$, but  $R(B_\lambda-\lambda)$ is not necessarily closed. Indeed, consider a normal operator $A\in B(H_i)$ such that  $\sigma(A)=[0,1]$ (for example, the multiplication operator $A:L^2([0,1])\to L^2([0,1])$ defined by $A(f)(x)=xf(x)$). Then $A\in\mathcal{C}(i)$. We show that for every $\lambda\in \sigma_p(A^\circ)$, $R(B_\lambda-\lambda)$ is not closed. By contradiction, suppose that there exists $\lambda\in\sigma_p(A^\circ)$ such that $R(B_\lambda-\lambda)$ is  closed. Then $B_\lambda-\lambda$ is a semi-Fredholm operator such that $\alpha(B_\lambda-\lambda)=0$. By \cite[Theorem 4.2.1]{MR0415345}, there exists $\epsilon>0$ such that if $|\gamma-\lambda|<\epsilon$ then $B_\lambda-\gamma\in \Phi_+(N(A^\circ-\lambda)^\perp)$ and $\alpha(B_\lambda-\gamma)=0$. This implies that $R(A^\circ-\gamma)=(\lambda-\gamma)N(A^\circ-\lambda)\oplus R(B_\lambda-\gamma)$ is closed and $\alpha(A^\circ-\gamma)=\alpha((\lambda-\gamma)I)+\alpha(B_\lambda-\gamma)=0$ for all $\gamma\in B(\lambda,\epsilon)$ with $\gamma\neq \lambda$. Therefore, $\lambda\in\iso\,\sigma_a(A^\circ)$. On the other hand, since $A$ is a normal operator it follows that $\sigma(A)=\sigma_a(A)=\sigma_a(A^\circ)$. Thus, $\lambda\in\iso\,\sigma(A)$, which is a contradiction, because $\sigma(A)=[0,1]$. Consequently, the equality (\ref{equation15}) is not true. This affects the proof of the main result in the paper \cite{MR2651679}. 

\end{rem}

\bibliographystyle{abbrv}

\bibliography{Manuscript}

\def\Dbar{\leavevmode\lower.6ex\hbox to 0pt{\hskip-.23ex \accent"16\hss}D}
  \def\dbar{\leavevmode\hbox to 0pt{\hskip.2ex \accent"16\hss}d}
\begin{thebibliography}{10}

\bibitem{MR1886113}
M.~Ahues, A.~Largillier, and B.~V. Limaye.
\newblock {\em Spectral computations for bounded operators}, volume~18 of {\em
  Applied Mathematics (Boca Raton)}.
\newblock Chapman \& Hall/CRC, Boca Raton, FL, 2001.

\bibitem{MR1047771}
A.~Aluthge.
\newblock On {$p$}-hyponormal operators for {$0<p<1$}.
\newblock {\em Integral Equations Operator Theory}, 13(3):307--315, 1990.

\bibitem{MR3624565}
A.~Ammar.
\newblock Some properties of the {W}olf and {W}eyl essential spectra of a
  sequence of linear operators {$\nu$}-convergent.
\newblock {\em Indag. Math. (N.S.)}, 28(2):424--435, 2017.

\bibitem{MR4075017}
A.~Ammar, A.~Bouchekoua, and A.~Jeribi.
\newblock Some approximation results in a non-{A}rchimedean {B}anach space.
\newblock {\em Funct. Anal. Approx. Comput.}, 12(1):33--50, 2020.

\bibitem{MR1212728}
L.~Burlando.
\newblock Spectral continuity in some {B}anach algebras.
\newblock {\em Rocky Mountain J. Math.}, 23(1):17--39, 1993.

\bibitem{MR0415345}
S.~R. Caradus, W.~E. Pfaffenberger, and B.~Yood.
\newblock {\em Calkin algebras and algebras of operators on {B}anach spaces}.
\newblock Marcel Dekker, Inc., New York, 1974.
\newblock Lecture Notes in Pure and Applied Mathematics, Vol. 9.

\bibitem{MR543882}
J.~B. Conway and B.~B. Morrel.
\newblock Operators that are points of spectral continuity.
\newblock {\em Integral Equations Operator Theory}, 2(2):174--198, 1979.

\bibitem{MR634170}
J.~B. Conway and B.~B. Morrel.
\newblock Operators that are points of spectral continuity. {II}.
\newblock {\em Integral Equations Operator Theory}, 4(4):459--503, 1981.

\bibitem{MR701023}
J.~B. Conway and B.~B. Morrel.
\newblock Operators that are points of spectral continuity. {III}.
\newblock {\em Integral Equations Operator Theory}, 6(3):319--344, 1983.

\bibitem{MR1457127}
S.~Djordjevi{\'c}.
\newblock The continuity of the essential approximative point spectrum.
\newblock {\em Facta Univ. Ser. Math. Inform.}, (10):97--104, 1995.

\bibitem{MR1814478}
S.~V. Djordjevi\'{c} and B.~P. Duggal.
\newblock Weyl's theorems and continuity of spectra in the class of
  {$p$}-hyponormal operators.
\newblock {\em Studia Math.}, 143(1):23--32, 2000.

\bibitem{MR1793814}
S.~V. Djordjevi{\'c} and Y.~M. Han.
\newblock Browder's theorems and spectral continuity.
\newblock {\em Glasg. Math. J.}, 42(3):479--486, 2000.

\bibitem{MR2651679}
B.~P. Duggal, I.~H. Jeon, and I.~H. Kim.
\newblock Continuity of the spectrum on a class of upper triangular operator
  matrices.
\newblock {\em J. Math. Anal. Appl.}, 370(2):584--587, 2010.

\bibitem{MR1785076}
I.~S. Hwang and W.~Y. Lee.
\newblock The spectrum is continuous on the set of {$p$}-hyponormal operators.
\newblock {\em Math. Z.}, 235(1):151--157, 2000.

\bibitem{MR433257}
G.~R. Luecke.
\newblock A note on spectral continuity and on spectral properties of
  essentially {$G_{1}$} operators.
\newblock {\em Pacific J. Math.}, 69(1):141--149, 1977.

\bibitem{MR1074574}
G.~J. Murphy.
\newblock {\em {$C^*$}-algebras and operator theory}.
\newblock Academic Press, Inc., Boston, MA, 1990.

\bibitem{MR0051441}
J.~D. Newburgh.
\newblock The variation of spectra.
\newblock {\em Duke Math. J.}, 18:165--176, 1951.

\bibitem{MR2894894}
S.~S{\'a}nchez-Perales and S.~V. Djordjevi{\'c}.
\newblock Continuity of spectra and compact perturbations.
\newblock {\em Bull. Korean Math. Soc.}, 48(6):1261--1270, 2011.

\bibitem{MR3388799}
S.~S\'{a}nchez-Perales and S.~V. Djordjevi\'{c}.
\newblock Spectral continuity using {$\nu$}-convergence.
\newblock {\em J. Math. Anal. Appl.}, 433(1):405--415, 2016.

\bibitem{MR3626681}
S.~S\'{a}nchez-Perales and S.~V. Djordjevi\'{c}.
\newblock Spectral continuity relative to invariant subspaces.
\newblock {\em Complex Anal. Oper. Theory}, 11(4):927--941, 2017.

\bibitem{sanchez4}
S.~S{\'a}nchez-Perales, S.~V. Djordjevi{\'c}, and S.~Palafox.
\newblock Some results about spectral continuity and compact perturbations.
\newblock {\em Filomat}, 2020.
\newblock Submitted.

\bibitem{MR4043866}
S.~S\'{a}nchez-Perales, S.~Palafox, and S.~V. Djordjevi\'{c}.
\newblock {$\alpha$}-{F}redholm operators relative to invariant subspaces.
\newblock {\em Oper. Matrices}, 13(4):921--936, 2019.

\bibitem{MR1861991}
M.~Schechter.
\newblock {\em Principles of functional analysis}, volume~36 of {\em Graduate
  Studies in Mathematics}.
\newblock American Mathematical Society, Providence, RI, second edition, 2002.

\bibitem{MR1670405}
A.~Uchiyama.
\newblock Berger-{S}haw's theorem for {$p$}-hyponormal operators.
\newblock {\em Integral Equations Operator Theory}, 33(2):221--230, 1999.

\end{thebibliography}

\end{document}